\documentclass[10pt]{article}  % 17/09/2017
\usepackage{amssymb}
\usepackage{amsmath,amsfonts,amsthm}

\newtheorem{theorem}{Theorem}[section]
\newtheorem{lemma}[theorem]{Lemma}
\newtheorem{proposition}[theorem]{Proposition}
\newtheorem{corollary}[theorem]{Corollary}

\newtheorem{definition}[theorem]{Definition}

\newcommand{\cqd}{\hfill $\Box$ \vspace{0.2cm}}

\usepackage{fullpage}

\begin{document}

\numberwithin{equation}{section}

\title{THE MEAN CURVATURE FLOW BY PARALLEL HYPERSURFACES}

\author{Hiuri Fellipe Santos dos Reis\footnote{Universidade de Bras\'ilia, Department of Mathematics, 70910-900, Bras\'ilia-DF, Brazil, hiuri.reis@ifg.edu.br  Partially supported  by CNPq Proc. 141275/2014-6, Ministry of Science and Technology, Brazil. } \qquad Keti Tenenblat\footnote{ Universidade de Bras\'{\i}lia,
 Department of Mathematics,
 70910-900, Bras\'{\i}lia-DF, Brazil, K.Tenenblat@mat.unb.br Partially supported by CNPq Proc. 312462/2014-0, Ministry of Science and Technology, Brazil.}
 }

\date{}

\maketitle{}

\begin{abstract}
It is shown that a hypersurface of a space form is the initial data for a solution to the mean curvature flow by parallel hypersurfaces if, and only if, it is isoparametric. By solving an ordinary 
differential equation, explicit solutions are given for all 
isoparametric hypersurfaces of space forms. In particular, for  
such hypersurfaces of the sphere, the exact collapsing time into a focal submanifold   is 
given in terms of its dimension, the principal curvatures and their multiplicities. 
\end{abstract}

\vspace{0.2cm} 
 
  %53C44 % local Riemannian geometry\\ 

\section{Introduction}

The Mean Curvature Flow (MCF) is a gradient-type flow for the 
volume functional. Under the MCF, a closed hypersurface in 
$\mathbb{R}^{n+1}$ locally evolves in the direction where the volume element decreases the fastest and eventually it becomes 
extinct. Along  the flow,  singularities may occur and one is interested in studying such singularities. There is an extensive literature on the subject starting with the early work in material science dated in the 1920s. We refer the reader to the excellent survey by Colding, Minicozzi and  Pedersen \cite{colding} and the references within. 

In recent years, self-similar solutions to the MCF have been studied. These are solutions given by the composition of isometries and homotheties. In the Euclidean space, the simplest examples of self-similar solutions are the spheres and the cylinders, which are self-contracting hypersurfaces. Other examples of self-contractile hypersurfaces in the Euclidean space can be found in \cite{angenent1, kleene} and translation hypersurfaces are found in \cite{clutterbuck,nguyen,wang,halldorsson1}. 
There are very few results on the mean curvature flow in non-Euclidean spaces. In \cite{hungerbuhler}, Hungerbühler and Smoczyk considered a particular case of self-similar solutions evolving by the MCF by a group of isometries of the ambient space, which are known as solitons, and presented several examples of these hypersurfaces on Riemannian manifolds. 
In \cite{liu}, Liu and Terng studied the MCF on isoparametric submanifolds with higher codimension of the Euclidean space and of the sphere, where they proved that the flow preserves the condition of being isoparametric and develops singularities in finite time, converging to a smooth submanifold of lower dimension. 

In this paper, we prove that any immersed hypersurface $M^n$ of a space form  evolves through the MCF  by parallel hypersurfaces if, and only if, $M$ is an isoparametric hypersurface. 
The MCF of such a hypersurface is obtained by solving an ordinary differential equation. We solve this  equation
for   all isoparametric hypersurfaces of space forms. 
In particular, we provide explicit solutions to the MCF of isoparametric hypersurfaces of the sphere and of the hyperbolic space,  including the exact  
collapsing time and the converging submanifold in terms of its dimension, the principal curvatures and their multiplicities. Our results are stated in the following section and  the  proofs are  given in Section \ref{proofs}.

%%%%%%%%%%%%%%%%%%%%%%%%%%%%%%%%%%%%

\section{Mean Curvature Flow of Isoparametric Hypersurfaces -  
Main results}\label{mainresults}

In what follows, $\mathbb{M}^{n+1}(\overline{\kappa})$ will be a space form of constant sectional curvature $\overline{\kappa} \in \left\{-1,0,1\right\}$, i. e.,  
$\mathbb{R}^{n+1}$ if $\overline{\kappa} = 0$, $\mathbb{S}^{n+1} \subset \mathbb{R}^{n+2}$ if $ \overline{\kappa} = 1$ 
and $\mathbb{H}^{n+1} \subset \mathbb{L}^{n+2}$ if $\overline{\kappa} = -1$, 
where $\mathbb{L}^{n+2}$ is the Lorentzian space. 
We consider   $F:M^n \rightarrow \mathbb{M}^{n+1}(\overline{\kappa})$ a hypersurface immersed in the space form $\mathbb{M}^{n+1}(\overline{\kappa})$,  with the induced metric 
$g(v,w) = \left\langle dF(v),dF(w)\right\rangle$, 
for all vector fields $v, \ w $ tangent to $M$. 
 If $F(M)$ is oriented and $N$ is a unit normal vector field, 
the \textit{second fundamental form} of $F(M)$ is given by  
 $h\left(v,w\right)= - \left\langle dN(v), dF(w) \right\rangle$.   
  Let $e_1, ...,e_n$ be orthonormal vector fields which are principal directions 
and let $\kappa_1, ..., \kappa_n$, be the \textit{principal curvatures} of $F(M)$
i.e., $g(e_\imath,e_\jmath) =  \delta_{\imath\jmath}$ and   $h(e_\imath,e_\jmath)=\kappa_\imath\delta_{\imath\jmath}$,
for  $1 \leq \imath, \jmath \leq n$. 
 We will denote the \textit{mean curvature} by $H=\sum_{\imath =1}^n\kappa_\imath$.  
 When the principal curvatures $\kappa_\imath$ of $F(M)$ do not depend on $x$, for all $\imath = 1, ..., n$, we say that $F(M)$ is an  \textit{isoparametric} hypersurface. From now on, we consider  connected
 hypersurfaces.

Let $F:M^n \rightarrow \mathbb{M}^{n+1}(\overline{\kappa})$  be an oriented hypersurface with a unit normal vector field $N$. A one parameter family of hypersurfaces $\widehat{F}:M^n\times I\rightarrow \mathbb{M}^{n+1}(\overline{\kappa})$, $I\subset \mathbb{R}$, is a solution to the \textit{ mean curvature flow}  (MCF) with initial condition $F$,  if 
\begin{equation}\label{(1)}
\left\{
\begin{array}{l}
\displaystyle{\frac{\partial}{\partial t} \widehat{F}(x,t) = \widehat{H}(x,t) \widehat{N}(x,t)},  \\
\widehat{F}(x,0) = F(x),
\end{array}
\right.
\end{equation}
where $\widehat{H}^t(.)=\widehat{H}(.,t)=\sum_{i=1}^n \widehat{k}_i^t$ is the mean curvature and $\widehat{N}^t(.)=\widehat{N}(.,t)$ is a  unit normal vector field  of $\widehat{F}^t(M)$. When $F$ is a minimal  hypersurface i.e. $H=0$,  then the family $\widehat{F}(t,x)=F(x)$
gives a trivial solution to the MCF.

In this paper, we consider a special type of solution to the MCF by  imposing that the hypersurfaces $\widehat{F}^t$ to be parallel. We first  introduce the following  notation 
\begin{equation}\label{cs} c (\xi) =\left\{ \begin{array}{ll}
1 ,& \ \text{if} \ \overline{\kappa} = 0,  \\
\cos(\xi) ,& \ \text{if} \ \overline{\kappa} = 1, \\ 
\cosh(\xi) ,& \ \text{if} \ \overline{\kappa} = -1, 
\end{array} \right.
\quad \text{and} \quad  s (\xi) =\left\{ \begin{array}{ll}
\xi ,& \ \text{if} \ \overline{\kappa} = 0  \\
\sin(\xi) ,& \ \text{if} \ \overline{\kappa} = 1, \\ 
\sinh(\xi) ,& \ \text{if} \ \overline{\kappa} = -1,
\end{array} \right. 
\end{equation}

\begin{definition}\rm{
Let $\widehat{F} : M^n\times I \rightarrow \mathbb{M}^{n+1}(\overline{\kappa})$ be a solution to the mean curvature flow in $\mathbb{M}^{n+1}(\overline{\kappa})$ with initial condition $F: M^n \rightarrow\mathbb{M}^{n+1}(\overline{\kappa})$. We say $\widehat{F}$ is a \textit{solution to the mean curvature flow by parallel hypersurfaces} if there is a function $\xi : I \rightarrow \mathbb{R}$, such that $\xi(0)=0$ and
\begin{equation}\label{eq3.1}
\widehat{F}^t(x) = c\big(\xi(t)\big) F(x)+s\big(\xi(t)\big)N(x),
\end{equation}
for all $t \in I$, where $c:\mathbb{R}\rightarrow \mathbb{R}$ and $s:\mathbb{R}\rightarrow \mathbb{R}$ are the functions defined in (\ref{cs}).} 
\end{definition}

We now state our main results. 
  
\begin{theorem}\label{teo3.1}
Let $F : M^n \rightarrow \mathbb{M}^{n+1}(\overline{\kappa})$ be a hypersurface in a space form $\mathbb{M}^{n+1}(\overline{\kappa})$. Then $F(M)$ is the initial data of a solution to the MCF by parallel hypersurfaces if, and only if, $F(M)$ is an isoparametric hypersurface. 
\end{theorem}
As a consequence of the proof of this theorem, given in Section 3,  one obtains the MCF of the isoparametric hypersurfaces of space forms by solving an ordinary differential equation. 
Namely, we prove the following 
\begin{corollary}\label{rem} Let  $F:M^n\rightarrow \mathbb{M}^{n+1}(\overline{\kappa})$ be an isoparametric hypersurface, with unit normal vector field $N$ and principal curvatures $\kappa_\imath$.  
Then the  solution to the MCF  with  initial data 
$F$ is  given by (\ref{eq3.1})
where $s$ and $c$ are the functions defined in  (\ref{cs}) and  $\xi(t)$ is the solution  of  
\[
 \xi'(t)=\sum_{\imath=1}^n\frac{\overline{\kappa}s\left(\xi\left(t\right)\right)
 +\kappa_\imath c\left(\xi\left(t\right)\right)}{c\left(\xi\left(t\right)\right)-\kappa_\imath s\left(\xi\left(t\right)\right)},   \qquad \xi(0)=0.
\] 
\end{corollary}

As an application, of Corollary \ref{rem}, we obtain explicitly the  MCF 
by parallel hypersurfaces  of the isoparametric hypersurfaces of $\mathbb{R}^{n+1}$  and of $\mathbb{H}^{n+1}$ in Propositions \ref{eucl}-\ref{cyl}. The MCFs for non minimal
hypersurface of $\mathbb{S}^{n+1}$ with $g$ distinct curvatures are given    
in Propositions \ref{umbilSn+1}-\ref{g6}.

For the sake of completeness we include the result for isoparametric hypersurfaces of the Euclidean space, without proof, since it is well known. 
\begin{proposition}\label{eucl}
Let $F:\mathbb{S}^m\times \mathbb{R}^{n-m}\rightarrow \mathbb{R}^{n+1}$, $m\neq 0$, be the immersion of a cylinder (or sphere if $m=n$) in the Euclidean space, with $m$ principal curvatures equal to $\kappa\neq 0$ and $n-m$ null principal curvatures. Then, the solution to the MCF with initial condition $F(M)$,  is 
\begin{equation} \label{hatFRn+1}
\widehat{F}^t(x)=F(x) + \frac{1-\sqrt{1-2m\kappa^2  t}}{\kappa}N(x),
\end{equation}
  for all $t\in (-\infty , t^*)$, where $t^*=\frac{1}{2m\kappa^2}$.  Moreover, if $m=n$ the solution collapses into a point at $t^*$. If $m\neq n$ then the solution collapses into an $(n-m)$-dimensional plane of $R^{n+1}$, at  $t^*$. 
\end{proposition}

\begin{proposition}\label{horo}
Let $F:\mathbb{R}^n\rightarrow \mathbb{H}^{n+1}\subset \mathbb{L}^{n+2}$ be the immersion of a horosphere in the hyperbolic space, with  unit normal vector field $N$ and all principal curvatures $\kappa = \pm 1$. Then, the solution to the MCF  with initial data $F$ is 
\begin{equation}\label{horosfera}
\widehat{F}^t(x) = \cosh(nt)F(x) + \kappa\sinh(nt)N(x),
\end{equation}
for all $t\in\mathbb{R}$. Moreover, $\widehat{F}^t(\mathbb{R}^n)$ is a horosphere for all $t\in \mathbb{R}$.
\end{proposition}

The totally umbilic hypersurfaces of the hyperbolic space, different from the horospheres, are treated in the following result. 

\begin{proposition}\label{sfcmu}
Let $F:M^n\rightarrow \mathbb{H}^{n+1}\subset \mathbb{L}^{n+2}$ be the immersion of a totally umbilic hypersurface in the hyperbolic space,  with unit normal vector field $N$ and all principal curvatures equal to $\kappa$ where  $\kappa\not\in \{0, \pm 1\}$. Then, the solution to the MCF with initial condition $F(M)$ is given by 
\begin{equation}\label{FCMUH}
\widehat{F}^t(x) = \frac{\kappa^2 e^{-nt}-\sqrt{1-\kappa^2+\kappa^2 e^{-2nt}}}{\kappa^2-1}F(x) + \frac{\kappa e^{-nt}-\kappa \sqrt{1-\kappa^2+\kappa^2e^{-2nt}}}{\kappa^2-1}N(x).
\end{equation}
\begin{enumerate}
\item If $0<|\kappa|<1$, then $\widehat{F}^t$ is defined for  $t\in\mathbb{R}$ and it converges to a totally geodesic $n$-dimensional manifold when $t \rightarrow + \infty$.
\item If $|\kappa|>1$ then $\widehat{F}^t$ is defined for $t\in (-\infty, t^*)$, where $t^*=\frac{1}{2n} \ln\left(\frac{\kappa^2}{\kappa^2-1}\right)$ and  it collapses to a point at $t^*$.
\end{enumerate}
\end{proposition}

For the hyperbolic cylinder, we have the following solution to the mean curvature flow. 

\begin{proposition}\label{cyl}
Let $F:\mathbb{S}^{m_1}\times \mathbb{H}^{m_2}\rightarrow \mathbb{H}^{n+1}\subset\mathbb{L}^{n+2}$ be the  immersion of a cylinder in  the hyperbolic space, with $m_1$ 
principal curvatures equal to $\kappa_1>1$ and $m_2$ principal curvatures equal to $\kappa_2$,  
such that $\kappa_1\kappa_2=1$. Then the solution to the MCF  with initial condition $F$,  is given by 
\begin{equation}\label{parallelhyperbolic}
\widehat{F}^t(x)=\cosh(\xi(t))F(x) + \sinh(\xi(t))N(x),
\end{equation}
 where
\begin{equation}\label{cscyl}
\cosh(2\xi(t))= \frac{a\ell(t)-2\sqrt{q(t)}}{a^2-4}\quad  \quad\sinh(2\xi(t)) =\frac{2\ell(t)-a\sqrt{q(t)}}{a^2-4}.
\end{equation} 
\begin{equation}\label{ell}
q(t)=\ell^2(t)-a^2+4; \quad 
\ell(t)=(a-b)e^{-2nt} + b, \quad a =\kappa_1+\kappa_2 \quad \mbox{and} \quad b=-\frac{m_1-m_2}{n}(\kappa_1-\kappa_2).
\end{equation}
 $\widehat{F}^t$ is 
defined for all $t\in (-\infty , t^*)$ where 
$t^*=\frac{1}{2n}\ln \frac{m_1\kappa_1^2+m_2}{m_1(\kappa_1^2-1)}$. 
and it  collapses into an $m_2$-dimensional focal  submanifold  at $t^*$.
\end{proposition}

We will now consider 
the  isoparametric hypersurfaces  of the sphere.    
Munzner \cite{munzner} showed that the number $g$ of distinct principal curvatures,  for an isoparametric hypersurface $M^n\subset \mathbb{S}^{n+1}$,   is restricted to be $1,\, 2,\, 3,\, 4$ or $6$.
Moreover, he showed that connected isoparametric hypersurfaces of the sphere can be extended to compact ones. 
Cartan \cite{cartan} classified these hypersurfaces when $g\leq 3$. 
If $g=1$ , then $M^n$ is a sphere obtained as the intersection 
of $\mathbb{S}^{n+1}$ with a hyperplane of $\mathbb{R}^{n+2}$.
 If $g=2$
 then $M^n$ must be the standard product of spheres $\mathbb{S}^l_{r_1}\times \mathbb{S}^{n-l}_{r_2} 
 \subset \mathbb{S}^{n+1}$, where $r_1^2+r_2^2=1$.    
When $g=3$, Cartan proved that there are only four distinct  
isoparametric hypersurfaces of $\mathbb{S}^{n+1}$  with three distinct 
principal curvatures. Their dimensions  are $n=3m$  where $m=1,\, 2,\, 4$ or $8$ and all the  principal curvatures have the same multiplicity $m$ . 
 The classification of the isoparametric hypersurfaces of the sphere 
with $g=4$ or $6$ is still not complete. However, 
when $g=4$, M\"unzner \cite{munzner} (see also Cecil-Chi-Jensen \cite{cecilchijensen})   
proved that the principal curvatures $\kappa_1,\kappa_2,\kappa_3,\kappa_4$ can be ordered so that their corresponding multiplicities satisfy $m_1=m_3$ and $m_2=m_4$. 
When $g=6$, Munzner \cite{munzner} showed that all the principal curvatures must have the same multiplicities $m$ and Abresch \cite{abresch} showed that $m=1$  or $m=2$. 

One can determine the principal curvatures of $M^n\subset \mathbb{S}^{n+1}$ 
up to a constant. In fact, for all $g=2,3,4,6$, let $a\in\mathbb{R}$      $-1<a<1$, and $a=\cos(gs)$, i.e., $0<s<\pi/g$. We consider  
\begin{equation}\label{kappagn}
\kappa_j=\cot\left(s+\frac{j-1}{g}\pi\right), \quad j=1,...,g. 
\end{equation}
The hypersurfaces for  different values of the constant $a$ are parallel
in $\mathbb{S}^{n+1}$.  

In the next propositions, we consider the solution $\widehat{F}^t$ to the  MCF, by parallel hypersurfaces, with initial data an isoparametric hypersurfaces $F$ of the sphere. Therefore, 
\begin{equation}\label{parallelsphere}
\widehat{F}^t\left(x\right)= \cos(\xi(t)) F\left(x\right)+\sin(\xi(t)) N\left(x\right). 
\end{equation}
The function $\xi(t)$ will be determined in the following results, according to the 
number $g$ of distinct principal curvatures.
  We start considering the case $g=1$, i.e.,  the umbilical hypersurfaces of  $\mathbb{S}^{n+1}$. 

\begin{proposition}\label{umbilSn+1}
Let $F:M^n\rightarrow \mathbb{S}^{n+1}\subset\mathbb{R}^{n+2}$ be the immersion of a totally umbilic hypersurface in $\mathbb{S}^{n+1}$,  with  unit normal vector field $N$ and all principal curvatures are equal to $\kappa \neq 0$. Then the solution to the MCF with $F$ as initial data, is given by 
\begin{equation}\label{FCMEE}
\widehat{F}^t(x)= \frac{\kappa^2e^{nt}+\sqrt{q(t)}}{\kappa^2+1}F(x) + \frac{\kappa e^{nt}-\kappa\sqrt{q(t)}}{\kappa^2+1}N(x),
\qquad \mbox{where} \qquad  
 q(t)=\kappa^2+1-\kappa^2e^{2nt}. 
\end{equation}
$\widehat{F}^t$ is defined for all $t \in (-\infty, t^*)$ where $t^*= \frac{1}{2n}\ln\left(\frac{\kappa^2+1}{\kappa^2}\right)$ 
and it collapses to a point at $t^*$.
\end{proposition}

\begin{proposition}\label{hopf}
Let $F:\mathbb{S}_{r_1}^l\times\mathbb{S}_{r_2}^{n-l}\rightarrow \mathbb{S}^{n+1}\subset{R}^{n+2}$ be an isoparametric hypersurface in $\mathbb{S}^{n+1}$, with two distinct principal curvatures $\kappa_1$ and  $\kappa_2$ with multiplicities $l$ and $n-l$ respectively. 
Then $\kappa_1\kappa_2=-1$ and assuming the immersion is not minimal, 
we may consider $\kappa_1>\sqrt{(n-l)/l}>1$. 
 The solution to the MCF with initial data  $F$,   is $\widehat{F}^t$ given by (\ref{parallelsphere}) 
 where 
\begin{equation}\label{csg2}
\cos(2\xi(t)) =  \frac{a\, q(t)+2\sqrt{a^2+4-q^2(t)}}{a^2+4}\qquad \qquad 
\sin(2\xi(t)) =  \frac{2q(t)-a\sqrt{a^2+4-q^2(t)}}{a^2+4}
\end{equation}
and 
\begin{equation}\label{notationg2}
a={\kappa_1+\kappa_2},\qquad b=-\,\frac{n-2l}{n}(\kappa_1-\kappa_2),\qquad   q(t)=(a+b)e^{2nt}-b.  
\end{equation}
$\widehat{F}^t$  is defined for all $t \in [0,t^*)$, where   
$t^*=\frac{1}{2n} \ln\left(\frac{l(\kappa_1^2+1)}{l(\kappa_1^2+1)-n}\right)$ and it collapses into an 
$(n-l)$-dimensional focal submanifold of $F$ at $t^*$.  
\end{proposition}

\begin{proposition}\label{g3}
Let $F:M^n\rightarrow \mathbb{S}^{n+1}\subset{R}^{n+2}$ be a non minimal  isoparametric 
hypersurface in $\mathbb{S}^{n+1}$,
 with unit normal vector field $N$ and three distinct principal curvatures. $\kappa_1,\kappa_2,\kappa_3$. Then all the principal curvatures have the same multiplicity $m$,  where  $m=1,\, 2,\, 4$ or  $8$, i.e. $n=3m$. Moreover, we may consider $\kappa_1>\sqrt{3}/3$ and   
 $\kappa_2$ and $\kappa_3$ given by (\ref{kappagn}). 
  The solution to the MCF  with initial data  $F$,  is  $\widehat{F}^t$  given by 
 (\ref{parallelsphere})
 where
\begin{equation}\label{csg3}
\begin{split}
\cos(3\xi(t))= \frac{a^2e^{9mt}+3\sqrt{q(t)}}{a^2+9}
,\qquad \qquad \qquad 
\sin(3\xi(t))= \frac{a(3e^{9mt}-\sqrt{q(t)})}{a^2+9},\\
a=\kappa_1+\kappa_2+\kappa_3=\frac{3\kappa_1(\kappa_1^2-3)}{3\kappa_1^2-1},\qquad \qquad q(t)= a^2+9-a^2e^{18mt}. 
\qquad \qquad \end{split}
\end{equation}
$\widehat{F}^t$ is defined for all $t \in [0,t^*)$, where $t^*={\displaystyle \frac{1}{18m}\ln\left(1+\frac{9}{a^2} \right)}$
%or equivalently where $\xi(t^*)=
%\mbox{arc}\cot \kappa_1$.  
and it collapses into a $2m$-dimensional focal submanifold of $F(M)$ at $t^*$. 
\end{proposition}

\begin{proposition}\label{g4}
Let $F:M^n\rightarrow \mathbb{S}^{n+1}\subset{R}^{n+2}$ be a non minimal  isoparametric 
hypersurface of $\mathbb{S}^{n+1}$,
 with unit normal vector field $N$ and four distinct principal curvatures $\kappa_j$, with multiplicities $m_j$, $j=1,2,3,4$. Then we may consider  
\begin{equation}\label{kappa4}
\kappa_1>1, \qquad \kappa_2=\frac{\kappa_1-1}{\kappa_1+1},\qquad
\kappa_3=\frac{-1}{\kappa_1},\qquad \kappa_4=\frac{-(\kappa_1+1)}{\kappa_1-1}, 
\end{equation}
where the multiplicities $m_j$ satisfy $m_1=m_3$ and  $m_2=m_4$, $n=2(m_1+m_2)$.    
The solution to the MCF with initial data  $F$,  is $\widehat{F}^t$ given by (\ref{parallelsphere}) 
 where 
\begin{equation}\label{csg4}
\cos(4\xi(t))=\frac{aq(t)+4\sqrt{a^2+16 -q^2(t)}}{a^2+16},
\qquad 
\sin(4\xi(t))= \frac{4q(t)-a\sqrt{a^2+16-q^2(t)}}{a^2+16},
\end{equation}  
\begin{equation}\label{ABg4}
a=\sum_{j=1}^4\kappa_j=\frac{\kappa_1^4-6\kappa_1^2+1}{\kappa1(\kappa^2-1)}\qquad, b=\frac{2(m_1-m_2)(\kappa_1^2+1)^2}{n\kappa_1(\kappa_1^2-1)} \qquad \mbox{and} \qquad 
q(t)=(a+b)e^{4nt}-b.
\end{equation}
Moreover, $\widehat{F}^t$ is defined for all $t \in [0,t^*)$, where $t^*={\displaystyle
\frac{1}{4n}\ln\left( \frac{b+\sqrt{a^2+16}}{a+b} \right)} $  and  
it collapses into $(m_1+2m_2)$-dimensional focal submanifold of $F(M)$.  
\end{proposition}

\begin{proposition}\label{g6}
Let $F:M^n\rightarrow \mathbb{S}^{n+1}\subset{R}^{n+2}$ be a non minimal isoparametric 
hypersurface in $\mathbb{S}^{n+1}$,
 with unit normal vector field $N$ and six distinct principal curvatures $\kappa_j$, $j=1,...,6$.  Then $n=6m$, where $m=1,2$, and we may consider $\kappa_1>\sqrt{3}$ and   
$\kappa_j$ given by (\ref{kappagn}), for $2\leq j\leq 6$.   
% \begin{equation}\label{kappa6}
%\kappa_1>\sqrt{3},\qquad \kappa_j=\frac{\kappa_1\cot((j-1)
%\frac{\pi{6})-1}{\kappa_1+\cot((j-1)\frac{\pi}{6})}, \quad j=2,...6.
%\end{equation}
 The  solution to the MCF  with initial data  $F$,  is $\widehat{F}^t$ given by (\ref{parallelsphere}) 
 where
\begin{equation}\label{csg6}
\cos(6\xi(t))=\frac{a^2e^{36mt}+6\sqrt{q(t)}}{a^2+36},
\qquad 
\sin(6\xi(t))= \frac{a\left( 6 e^{36mt}-\sqrt{q(t)}\right)}{a^2+36},
\end{equation}
where  
\begin{equation}\label{ABg6}
a=\sum_{j=1}^6 \kappa_j=\frac{\kappa_1^6-15\kappa_1^4+15\kappa_1^2-1}{k_1(\kappa_1^2-3)(3\kappa_1^2-1)} \qquad \mbox{and} \qquad 
q(t)=a^2+36-a^2e^{72mt}.
\end{equation}
which is defined for all $t \in [0,t^*)$, where $t^*=\displaystyle{ 
\frac{1}{72m}\ln \left(1+\frac{36}{a^2}\right)}$. 
Moreover, the solution collapses into a $5m$-dimensional focal submanifold  of $F(M)$ at $t^*$. 
\end{proposition}

% VER FIGURAS

\section{Proof of the main results}\label{proofs}
In order to prove our main results, 
we first state some well known properties of parallel hypersurfaces.
In fact, our next lemma can be easily proved as a consequence of 
(\ref{eq3.1}) and the fact that the functions $c(\xi)$ and $s(\xi)$, defined by (\ref{cs}),  satisfy the following properties
\begin{equation}\label{properties}
c'(\xi)=-\overline{\kappa}s(\xi), \qquad 
s'(\xi) = c(\xi), \qquad 
c^2(\xi)+\overline{\kappa}s^2(\xi) =1. 
\end{equation}
\begin{lemma}\label{lema1}
Let $F : M^n\rightarrow \mathbb{M}^{n+1}(\overline{\kappa})$  be an oriented hypersurface with unit normal vector field $N$ and  let $\widetilde{F}^\xi:M^n\rightarrow \mathbb{M}^{n+1}(\overline{\kappa})$ be a parallel hypersurface, given by $\widetilde{F}^\xi(x)= c(\xi)F(x)+s(\xi)N(x)$. 
 Then, the  unit  vector field $\widetilde{N}^\xi$ normal to   $\widetilde{F}^\xi$  
   and the corresponding principal curvatures $  \widetilde{\kappa}_\imath^\xi$ are given by 
\begin{equation}\label{NPki}
\widetilde{N}^\xi (x)=-\overline{\kappa} s(\xi) F(x) + c(\xi) N(x), \qquad \qquad 
\widetilde{\kappa}_\imath^\xi(x) = \frac{\overline{\kappa} s\big(\xi\big) + \kappa_\imath(x) c\big(\xi\big)}{c\big(\xi\big) - \kappa_\imath(x) s\big(\xi\big)}, 
\end{equation}
where 
$\kappa_{\imath}(x)$ is the $\imath$-th principal curvature of $F(M)$ on $x\in M$, for all $\imath = 1,...,n$. 
Moreover, if $\{e_1,...,e_n\} \subset T_xM$ is an orthonormal basis of eigenvectors of the second fundamental form of $F(M)$, then 
\begin{equation}\label{gij}
\widetilde{g}^\xi_x(e_\imath, e_\jmath) = \big[c(\xi)-\kappa_\imath(x)s(\xi)\big]^2\delta_{\imath \jmath}.
\end{equation} 
\end{lemma}

\noindent\textbf{ Proof of Theorem \ref{teo3.1}.} Suppose that $\widehat{F}^t$ is given by (\ref{eq3.1}). Then,  as a consequence of  (\ref{properties}), 
\begin{eqnarray}\label{deriv}
\frac{\partial }{\partial t}\widehat{F}(x,t) 
= \xi'(t)\big(-\overline{\kappa} s(\xi(t))F(x) + c(\xi(t))N(x) \big),    
\end{eqnarray}
If $\widehat{F}^t$   is a solution to MCF, then   (\ref{(1)}) reduces to   
$$\xi'(t)\big(-\overline{\kappa} s(\xi(t)) F(x)+c(\xi(t)) N(x)\big)=\widehat{H}^t(x)\widehat{N}^t(x)$$
and  from (\ref{NPki}) we conclude  that 
\begin{equation}\label{Hctt}
\xi'(t) = \widehat{H}^t(x). 
\end{equation}
Therefore, for each $t$ fixed, the mean curvature $\widehat{H}^t$ is constant. Then, it follows from a classical result of Cartan \cite{cartan} (see also Theorem 3.6 in \cite{cecil}) that $F(M)$ is an isoparametric hypersurface.

Reciprocally, if $F(M)$ is an isoparametric hypersurface,  let $\xi(t)$ be the unique solution of  the ordinary differential equation
\begin{eqnarray}\label{edop}
\xi'(t)=\sum_{\imath=1}^n\widehat{\kappa}^t_\imath=
\sum_{\imath=1}^n\frac{\overline{\kappa}s\left(\xi\left(t\right)\right)+\kappa_\imath c\left(\xi\left(t\right)\right)}{c\left(\xi\left(t\right)\right)-\kappa_\imath s\left(\xi\left(t\right)\right)},
\end{eqnarray}
such that  $\xi(0)=0$, where $\kappa_\jmath$ is the $\jmath$-th principal curvature of $F(M)$. 

Considering $\widehat{F}^t(x)$ given by (\ref{eq3.1}), then (\ref{deriv}) 
holds. Moreover,  
it follows from  Lemma \ref{lema1}, that $$\widehat{N}^t(x) =-\overline{\kappa}s(\xi(t)) F(x) + c(\xi(t))N(x)$$
and 
\begin{eqnarray*}
\widehat{H}^t(x) = \sum_{\imath=1}^n\widehat{\kappa}^t_\imath=\sum_{\imath=1}^n\frac{\overline{\kappa}s\left(\xi\left(t\right)\right)+\kappa_\imath c\left(\left(t\right)\right)}{c\left(\left(t\right)\right)-\kappa_\imath s\left(\left(t\right)\right)} = \xi'(t).
\end{eqnarray*}
Therefore, 
\begin{eqnarray*}
\frac{\partial}{\partial t}\widehat{F}^t(x) &=& \widehat{H}^t(x) \widehat{N}^t(x). 
\end{eqnarray*}
Moreover since  $\xi(0)=0$, we have $\widehat{F}^0(x)=F(x)$, i.e., 
$\widehat{F}^t$ is a solution to the MCF.
\cqd
 
The proof of  Corollary \ref{rem} is an immediate  consequence of the proof given above. 

From now on, we will  apply  Corollary \ref{rem} to obtain the mean curvature flow whose  initial data is a  non minimal   
isoparametric hypersurface of a space form, i.e. we will prove Propositions 
\ref{eucl}-\ref{g6}.

%\noindent\textbf{Proof.} From Corollary \ref{rem}, the solution to the  %mean curvature flow, by parallel hypersurfaces, is given by %$\widehat{F}^t(x)= F(x)+\xi(t)N(x),$ where $\xi(t)$ satisfies 
%$\xi'(t) = m\kappa/(1-\kappa \xi(t))$ with $\xi(0)=0$.  
%Integrating this equation, we obtain (\ref{hatFRn+1}) 
%for all  $t\in \left(\infty, 1/2m\kappa^2\right)$. Moreover, it follows %from (\ref{NP}) that $\widehat{N}=N$. When $t^*=1/2m\kappa^2$,  
%$\widehat{F}^t(x)$ collapses into $\widehat{F}^{t^*(x)
%=F(x)+\frac{1}%{\kappa}N(x)$.
%If $m=n$, since $dF_x=-dN_x/\kappa$ $\widehat{F}^*(x)$ reduces to  a %point.  
%since $dF_x = -dN_x/\kappa$, we have 
%$d\widehat{F}^{\frac{1}{2n\kappa^2}}_x= 0$.  Therefore, %$\widehat{F}^{\frac{1}{2n\kappa^2}}(x)$ is a point. 
%If $n\neq m$, then for all 
% $V_1\in T\mathbb{S}^m$ and $V_2\in T\mathbb{R}^{n-m}$, we have  
% $d\widehat{F}^{t^*}_x(V_1)=0,\quad 
%d\widehat{F}^{t^*}_x(V_2) = dF_x(V_2)$ and 
%$d\widehat{N}_x(V_2) =dN_x(V_2)=0$. 
%Thus, $\widehat{F}^{t^*}$ is an $(n-m)$-dimensional plane in % %    $\mathbb{R}^{n+1}$.
%\cqd

\noindent\textbf{Proof of Proposition \ref{horo}.} From Corollay \ref{rem} and the fact that $\kappa=\pm 1$, we have $\widehat{F}^t(x)= \cosh(\xi(t))F(x)+\sinh(\xi(t))N(x)$, where $\xi'(t)=\kappa n$. 
Integrating this equation, since $\xi(0)=0$, we have  $\xi(t) =\kappa nt,$ for all $t\in \mathbb{R}$. Therefore, (\ref{horosfera}) follows.
Moreover, from (\ref{NPki}), we have 
$\widehat{\kappa}^t_\imath = \kappa$, $\forall \imath$, which concludes the proof.

\cqd

\noindent\textbf{Proof of Proposition \ref{sfcmu}.} From Corollary \ref{rem}, the solution to the MCF by parallel hypersurfaces is given by  $$\widehat{F}^t(x)= \cosh(\xi(t))F(x)+\sinh(\xi(t))N(x),\quad  \mbox{where} \quad  \xi'(t)=-n\frac{\sinh(\xi(t))-\kappa\cosh(\xi(t))}{\cosh(\xi(t))-\kappa\sinh(\xi(t))}.$$
Integrating the equation for $\xi$, with $\xi(0)=0$, we have 
$\sinh(\xi(t))=\kappa\cosh(\xi(t))-\kappa e^{-nt}$.  
The square of this equation reduces to  
$$(\kappa^2-1)\cosh^2(\xi(t))-2\kappa^2e^{-nt}\cosh(\xi(t))+\kappa^2e^{-2nt}+1=0.$$ 
Therefore, 
$$\cosh\left(\xi\left(t\right)\right) = \frac{\kappa^2 e^{-nt}-\sqrt{q(t)}}{\kappa^2-1},\quad \sinh\left(\xi\left(t\right)\right) =\frac{\kappa e^{-nt}-\kappa \sqrt{q(t)}}{\kappa^2-1}, \quad \mbox{  where  }  q(t)=1-\kappa^2+\kappa^2 e^{-2nt},$$
%$$\cosh\left(\xi\left(t\right)\right) = \frac{\kappa^2 
%e^{-nt}-\sqrt{1-\kappa^2+\kappa^2 e^{-2nt}}}{\kappa^2-1} \quad %\text{and}\quad \sinh\left(\xi\left(t\right)\right) =\frac{\kappa 
%e^{-nt}-\kappa \sqrt{1-\kappa^2+\kappa^2e^{-2nt}}}{\kappa^2-1},$$
which proves (\ref{FCMUH}). 
 Let $\{E_1,...,E_n\}$ be a local orthonormal frame of principal directions on $F(M)$. It follows from Lemma \ref{lema1}, that the first fundamental form and the principal curvatures of $\widehat{F}^t$  are given respectively  by $\widehat{g}^t(E_\imath,E_\jmath)=q(t)\delta_{\imath\jmath}$,
and $\widehat{\kappa}^t= \kappa e^{-nt}/\sqrt{(q(t)}$.

If $0<|\kappa|<1$, then $q(t)>0$ and   
 the solution is defined  for all $t\in\mathbb{R}$.
  Moreover,since $\lim_{t\rightarrow +\infty}\widehat{\kappa}^t =0$, 
 $\widehat{F}^t(\mathbb{H}^n)$ converges to a 
  totally geodesic submanifold.     
  If $|\kappa|>1$, then $ q(t^*)=0$ where $t^*=\frac{1}{2n} \ln\left(\frac{\kappa^2}{\kappa^2-1}\right)$. Hence 
$\widehat{F}^t(\mathbb{S}^n)$  is defined for $t\in (-\infty, t^*)$ and it collapses to a point at $t^*$.
  
\cqd

\noindent\textbf{Proof of Proposition \ref{cyl}.} It follows from Corollary \ref{rem}  that the solution $\widehat{F}^t$ is given by  (\ref{parallelhyperbolic}) 
with   $\xi(t)$ satisfying  
$$\xi'(t)=-n\frac{2\sinh(2\xi(t))-a\cosh(2\xi(t))+b}{2\cosh(2\xi(t))-a\sinh(2\xi(t))},$$
where $a$ and $b$ are given by (\ref{ell}).
Integrating, this equation whith $\xi(0)=0$, we get
$2\sinh(2\xi(t))=a\cosh(2\xi(t))-b+(b-a)e^{-2nt}$. 
The square of this equation   reduces to 
\begin{eqnarray*}
(a^2-4)\cosh^2(2\xi(t))-2a\ell(t)\cosh(2\xi(t)) + \ell^2(t)+4=0,
\end{eqnarray*}
where $\ell(t)$ is given by (\ref{ell}). 
Therefore, we obtain $\cosh(2\xi(t))$ and $\sinh(2\xi(t))$ as in (\ref{cscyl}). Without loss of generality, we are considering $\kappa_1>\kappa_2>0$. 

Moreover, $\widehat{F}^t$ is defined for all $t$, for which  
 $q(t)=\ell^2(t)-a^2+4 > 0$. 
Since $-a^2 +4= -(\kappa_1-\kappa_2)^2$ and  $b-a=-{2}(m_1\kappa_1+m_2\kappa_2)/n $, we conclude that $t<t^*$, where $t^*=1/n \ln\frac{m_1\kappa_1+m_2\kappa_2)}{m_1(\kappa_1-\kappa_2)} $.

Let $\{E_1,...,E_{m_1},E_{m_1+1},...,E_{m_1+m_2}=E_n\}$ be an orthonormal frame of principal directions such that $E_1,...,E_{m_1}$ are tangent to $\mathbb{S}^{m_1}$ and $E_{m_1+1},...,E_{n}$ are tangent to $\mathbb{H}^{m_2}$. It follows from (\ref{gij}) that $\widehat{g}^{t^*}(E_\imath, E_\jmath)=0$ for $1\leq \imath\leq m_1$ and  $\widehat{g}^{t^*}(E_\imath, E_\jmath)=\kappa_2(\kappa_1-\kappa_2)\delta_{\imath\jmath}$ for  
 $m_1+1\leq \imath,\jmath\leq n$.  
Thus, the solution $\widehat{F}^t$ collapses into an $m_2$-dimensional focal submanifold at $t^*$, 
since $\coth (\xi(t^*))=\kappa_1$. 

\cqd

We will now prove Propositions \ref{umbilSn+1}-\ref{g6}  that consider the MCF of isoparametric hypersurfaces of the sphere.   

\noindent\textbf{Proof of Proposition \ref{umbilSn+1}.} It follows from Corollary \ref{rem}, that the solution to the MCF is given by  
(\ref{parallelsphere}),  where $\xi(t)$ satisfies     
$$\xi'\left(t\right) = n\, \frac{\sin\left(\xi\left(t\right)\right)+\kappa\cos\left(\xi\left(t\right)\right)}{\cos\left(\xi\left(t\right)\right)-\kappa\sin\left(\xi\left(t\right)\right)},\qquad \mbox{with}\qquad \xi(0)=0.$$ 
Integrating this equation, we have
$\sin\left(\xi\left(t\right)\right)=
\kappa( e^{nt}-\cos\left(\xi\left(t\right)\right))$. 
The square of this equation reduces to 
$$\left(\kappa^2+1\right)\cos^2\left(\xi\left(t\right)\right) -2 \kappa^2e^{nt}\cos\left(\xi\left(t\right)\right) +\kappa^2e^{2nt}-1=0.$$
Therefore, we obtain (\ref{FCMEE}) and the solution is defined for all
$t\in (-\infty, t^*)$ where 
 $t^*= \frac{1}{2n}\ln\left(\frac{\kappa^2+1}{\kappa^2}\right)$. 
Let $\left\{ E_1, ..., E_n\right\}$ be a local orthonormal frame  of principal directions of $F(M)$. From Lemma \ref{lema1}, we have
$\widehat{g}^t\left(E_\imath,E_\jmath\right)= 
=q^2(t)\delta_{\imath\jmath}$.  
Since $q(t^*)=0$, we have   $\widehat{g}^{t^*}(E_\imath,E_\jmath) =0$,  for all $\imath,\jmath$. Therefore, $\widehat{F}^{t}(M)$ collapses into a point at $t^*$.

\cqd

\noindent\textbf{Proof of Proposition \ref{hopf}.} 
Since we are considering that $F$ is not a minimal hypersurface, we may assume that  the mean curvature $H>0$.   Since (\ref{kappagn}) implies that  $\kappa_2=-1/\kappa_1$, hence,  
without loss of generality,  we may consider  $\kappa_1>\sqrt{(n-l)/l\,}$ and $n-l>l$ i.e. $\kappa_1>1$.      
It follows from  Corollary \ref{rem} and  a straightforward computation that the solution of the mean curvature flow is given by
 (\ref{parallelsphere}), 
where $\xi(t)$ must satisfy 
\[
\xi'(t)= n \,\frac{2\sin(2\xi(t)) +a \cos(2\xi(t))+b}{2\cos(2\xi(t)) -a\sin(2\xi(t))}, \qquad \xi(0)=0, 
\]
where $a$ and $b$ are given by (\ref{notationg2}). 
Integrating,  we get 
\begin{equation}\label{sing2}
2\sin(2\xi(t)) + a\cos(2\xi(t))+b = (a+b)e^{2nt}.
\end{equation} 
This equation implies that $\cos(2\xi(t)$ must satisfy the algebraic equation  
$$(a^2+4)\cos^2(2\xi(t))-2aq(t)\cos(2\xi(t)) + q^2(t)-4=0,$$
where $q(t)$ is given by (\ref{notationg2}). 
Solving for $\cos(2\xi(t))$ with $\xi(0)=0$ and using (\ref{sing2}), we get  (\ref{csg2}).   
 
Let $t^*$ be such that $a^2+4-q^2(t^*)=0$. Then, since $a+b=2H/n>0$, we get  
\[
t^*=\frac{1}{2n}\ln\left(\frac{b+\sqrt{a^2+4}}{a+b} \right)=
\frac{1}{2n}
\ln\left(\frac{l(\kappa_1^2+1)}{l(\kappa_1^2+1)-n}  \right),  
\]
 and it follows from (\ref{csg2}) that $\cot (2\xi(t^*))=a/2$ and 
$\cot(\xi(t^*))=\kappa_1$. 
Let $\{E_1,...,E_l\}$ and $\{E_{l+1},...,E_n\}$ be an orthonormal frame of principal vector fields  
corresponding to $\kappa_1$ and $\kappa_2$ respectively. 
It follows from Lemma \ref{lema1} that $\widehat{g}^{t^*}(E_i,E_j)=(\sin(\xi(t^*))^2(\cot(\xi(t^*))-\kappa_1)^2\delta_{ij}=0$, for 
$1\leq i,j\leq l$ and $\widehat{g}^{t^*}(E_i,E_j)= \sin(\xi(t^*))^2(\cot(\xi(t^*))+1/\kappa_1)^2\delta_{ij}\neq 0$, for 
$l+1\leq i,j\leq n$. Therefore, Hence $\widehat{F}^t$ is defined for $t\in [0,t^*)$ and when $t$ tends to $t^*$, the mean curvature flow collapses into an $(n-l)$-dimensional focal submanifold
of $F(M)$. 

\cqd

%g=3 n=3m onde m=1,2,4 ou 8

\noindent\textbf{Proof of Proposition \ref{g3}.} 
When $g=3$ the three distinct principal curvatures are determined by 
(\ref{kappagn}). Since $H>0$ we may consider   $
 \kappa_1>\sqrt{3}/{3}$. Then 
$ \kappa_2= (\kappa_1\sqrt{3}-3)/(3\kappa_1+\sqrt{3})$ and 
  $\kappa_3= (\kappa_1\sqrt{3}+3)/(-3\kappa_1+\sqrt{3})$. 
Moreover, Cartan showed that all the principal curvatures have the same multiplicity 
$m$, where  $m=1,\, 2,\, 4$ or 
 $8$, and hence  $n=3m$. 
It follows from  Corollary \ref{rem} and  a straightforward computation that the solution of the mean curvature flow is given by
(\ref{parallelsphere}), 
 where $\xi(t)$ must satisfy 
\[
\xi'(t)= 3m\frac{3\sin(3\xi(t))+a\cos(3\xi(t))}{3\cos(3\xi(t))+a\sin(3\xi(t))}, 
\]
where $a=\sum^3_{\jmath=1}\kappa_\jmath$. 
Integrating this equation, and using the fact that $\xi(0)=0$,  we get 
\begin{equation}\label{cg3}
3\sin(3\xi(t))+a\cos(3\xi(t))=ae^{9mt}.
\end{equation}
Considering the square of this equation, we obtain
\[
(a^2+9)\cos^2(3\xi(t))-2a^2e^{9mt}\cos(3\xi(t))+a^2e^{18mt}-9=0.
\]
Solving this equation for $\cos(3\xi(t))$ with $\xi(0)=0$ and using  (\ref{cg3}) we obtain $\cos(3\xi(t))$ 
and $\sin(3\xi(t))$ given by (\ref{csg3}).  
Let $t^*$ be such that $q(t^*)=0$. Then 
\[
t^*= \frac{1}{18m}\ln\left(1+\frac{9}{a^2} \right)=
\frac{1}{18m}\ln\left(\frac{(\kappa_1^2+2)^3}{\kappa_1^2(\kappa_1^2-3)^2} \right) 
\]
 and it follows from (\ref{csg3}) that   $\cot (3\xi(t^*))=a/3$.
Moreover, $\cot(\xi(t^*))=\kappa_1$. 
%In fact, let $t_1$ be such 
%$\cot(\xi(t_1))=\kappa_1$. Then a simple calculation shows that 
%\[
%\cot(\xi(3t_1))=\frac{\kappa_1(\kappa_1^2-3)}{3\kappa_1^2-1}= \frac{a}{3}.
%\]
%Hence, $t_1=t^*$ and   $\cot(\xi(t^*))=\kappa_1$. 
When $t$ tends to $t^*$ 
the hypersurface $\hat{F}^t$ collapses into a $2m$-dimensional focal submanifold of 
$F(M)$. In fact, considering $E_{j1},...E_{jm}$, $j=1,2,3$ 
an orthonormal vector field of principal direction corresponding to the 
principal curvature $\kappa_j$, we conclude that 
$\hat{g}^{t^*}(E_{1\ell},E_{1r})=\sin(\xi(t^*)) (\cot (\xi(t^*))-\kappa_1)=0$ for all $1\leq \ell,r\leq m$, while $\hat{g}^{t^*}(E_{j\ell},E_{j\ell})\neq 0$ for $j\neq 1$. 

 \cqd

\noindent\textbf{Proof of Proposition \ref{g4}.}
 Since $g=4$, it follows from (\ref{kappagn})  that we may consider the principal curvature  given by (\ref{kappa4}). Moreover, M\"unzner proved that the corresponding multiplicities satisfy $m_1=m_3$ and $m_2=m_4$, hence 
 $n=2(m_1+m_2)$. It follows from Corollary \ref{rem} and a straightforward computation that the solution of the MCF $\widehat{F}^t$ is given by (\ref{parallelsphere}), where $\xi(t)$  must satisfy 
 $$
 \xi '(t) = n\frac{4\sin(4\xi(t))+a\cos(4\xi(t))+b}{4cos(4\xi(t))-a\sin(4\xi(t))}, \ \ \ \xi(0)=0,
 $$
  where $a$ and $b$ are the constants given by (\ref{ABg4}). Integrating, we get 
\begin{equation}\label{eq}
4\sin(4\xi(t))+a\cos(4\xi(t)) +b=(a+b)e^{4nt}.
\end{equation}  
It follows from this equation that $\cos(4\xi(t))$ must satisfy the algebraic equation
$$
(a^2+16)\cos^2(4\xi(t))-2a\big((a+b)e^{4nt}-b\big)\cos(4\xi(t))+\big((a+b)e^{4nt}-b\big)^2-16=0.
$$
Solving this equation for $\cos(4\xi)$ and considering $\xi(0)=0$,  the expression $q(t)$ defined by (\ref{ABg4}) and (\ref{eq}), we conclude that $\cos(4\xi(t))$ and $\sin(4\xi(t))$ are given by (\ref{csg4}). 

Let $t^*$ be such that $a^2+16-q^2(t^*)=0$. Then 
$$
t^*=\frac{1}{4n}\ln\bigg(\frac{b+\sqrt{a^2+16}}{a+b}\bigg)=\frac{1}{4n}\ln\Bigg(\frac{m_1(\kappa_1^2+1)^2}{m_1(\kappa_1^2+1)^2-2n\kappa_1^2}\Bigg).
$$
It follows from (\ref{csg4}) that $\cot(4\xi(t^*))=\frac{a}{4}$ and 
$ \cot(\xi(t^*)=\kappa_1$. 

 When t tends to $t^*$, the hypersurface collapses into an $(n-m_1)$-dimensional submanifold of $\mathbb{S}^{n+1}$, which is a focal submanifold of $F(M)$. In fact, if  we consider $\{ E_1,...,E_{m_1} \}$, 
  $\{ E_{m_1+1},...,E_{m_1+m_2} \}$, $\{ E_{m_1+m_2+1},...,E_{2m_1+m_2}\}$ and $\{  E_{2m_1+m_2+1},...,E_{n}\}$ 
 an orthonormal frame of principal directions of $F{M}$, corresponding to $\kappa_1, \ \kappa_2, \ \kappa_3, \ \kappa_4$, respectively, then 
 $\lim_{t\rightarrow t^*}\widehat{g}^{t}(E_{\imath}, E_{\jmath}) = 0$
 for all $\imath,\jmath $ such that $1\leq \imath,\jmath\leq m_1$ 
 and $\widehat{g}^{t^*}(E_\ell,E_r)=0$ for $\ell\neq r$ and $\ell>m_1$ or $r>m_1$ while , $\widehat{g}^{t^*}(E_\ell,E_\ell) \neq 0$, for $m_1+1\leq \ell\leq n$.   
 
\cqd

\noindent\textbf{Proof of Proposition \ref{g6}.} 
When $g=6$ it follows from results of Munzner and Abresch that all the principal curvatures have the same multiplicity $m$ and $m=1,2$. Moreover, since $F$ is not a minimal hypersurface we may consider 
$\kappa_1>\sqrt{3}$ and $\kappa_j$ are given by 
 (\ref{kappagn}) for $2\leq j\leq 6$. 
It follows from  Corollary \ref{rem} and  a straightforward computation that the solution is given by (\ref{parallelsphere}), 
 where $\xi(t)$ must satisfy 
\[
\xi'(t)=\frac{6m(a\cos(6\xi(t))+6\sin(6\xi(t)))}{6\cos(6\xi(t))-a\sin(6\xi(t))}\qquad \xi'(0)=0.
\]
Integrating this equation we get 
%\begin{equation}\label{cg6}
$ a\cos(6\xi(t))+6\sin(6\xi(t))=ae^{36mt}$ .
%\end{equation}
It follows from this equation that $\cos(6\xi(t))$ must satisfy the algebraic equation
\[
(a^2+36)\cos^2(6\xi(t))-2a^2e^{36mt}\cos(6\xi(t))+
a^2e^{72mt}-36=0.
\]
Solving this equation for $\cos(6\xi(t))$ and considering  $\xi(0)=0$,  
we conclude that $\cos(6\xi(t))$ and $\sin(6\xi(t))$ are given by 
 (\ref{csg6}). 
Let  $t^*$ be such that $q(t^*)=0$. Then $t^*=\displaystyle{ 
\frac{1}{72m}\ln \left(1+\frac{36}{a^2}\right)}$ and it follows from 
(\ref{csg6}) that $\cot(6\xi(t^*))=a/6$.  Moreover,   $\cot(\xi(t^*))=\kappa_1$. 
%In fact, let 
%$t_1$ be such that $\cot(\xi(t_1))=\kappa_1$. Then a simple calculation 
%shows that $\cot(6\xi(t_1))=a/(2b)$ and hence $t_1=t^*$. Therefore, 
%$\cot(t^*)=\kappa_1$.

Then $\widehat{F}^t$ is defined for $t\in [0,t^*)$ and  when $t$ tends to $t^*$ the hypersurface it 
collapses into a $5m$-dimensional submanifold of $\mathbb{S}^{6m+1}$
that is a focal submanifold of $F(M)$. 
 In fact, $E_{j1},E_{jm}$, $j=1,...,6$ an orthonormal  vector field of principal directions corresponding to the principal vurvature $\kappa_j$, we conclude that $\hat{g}^{t^*}(E_{1\ell}, E_{1r})=0$, for all $1\leq \ell,r\leq m$, while $\hat{g}^{t^*}(E_{j\ell},E_{j\ell})\neq 0$, for $j\neq 1$. 
  
\cqd

\end{document}